\documentclass[11pt, leqno]{article}
\usepackage{pdfsync}
\usepackage{array} 
\usepackage{amsfonts}
\usepackage{amsmath}
\usepackage{amssymb}
\usepackage{graphicx}
\usepackage{color}
\usepackage{bm}
\usepackage[utf8]{inputenc}
\usepackage[T1]{fontenc}
\newtheorem{theorem}{Theorem}[section]
\newtheorem{corollary}{Corollary}[theorem]

\newtheorem{proposition}[theorem]{Proposition}

\newcommand{\RR}{\mathbb{R}}
\newcommand{\re}{\mathbb{R}}

\newcommand{\ren}{\mathbb{R}^N}
\newcommand{\rd}{{\rm d}}
\newcommand{\ou}{{\widehat u}}
\newcommand{\ve}{{\varepsilon}}

\def\qed{\unskip\kern 6pt \penalty 500
\raise -2pt\hbox{\vrule \vbox to8pt{\hrule width 6pt \vfill\hrule}\vrule}\par}

\numberwithin{equation}{section}

\newcommand{\hiddensubsection}[1]
{    \stepcounter{subsection}
      \subsection*{\arabic{section}.\arabic{subsection}\hspace{1em}{#1}}
      }
\begin{document}

\title{\textbf{Asymptotic behaviour  methods for the Heat \\Equation.
Convergence to the Gaussian}}

\author{{\Large Juan Luis V\'azquez}\\ [4pt]
Universidad Aut\'{o}noma de Madrid\\
} 
\date{\today} 

\date{2018}

\maketitle

\begin{abstract}
\noindent In this expository work we discuss the  asymptotic behaviour of the solutions of the classical heat equation posed in the whole Euclidean space.
 After an introductory review of the main facts on the existence and properties of solutions, we proceed with the proofs of convergence to the Gaussian fundamental solution, a result that holds for all integrable solutions, and represents in the PDE setting the Central Limit Theorem of probability.  We present several methods of proof: first, the scaling method. Then several versions of the representation method. This is followed by the functional analysis approach that leads to the famous related equations, Fokker-Planck and Ornstein-Uhlenbeck. The analysis of this  connection is also given in rather complete form here.   Finally, we present the Boltzmann entropy method, coming from kinetic equations.\\
 \qquad The different methods are interesting because of  the possible extension to prove the asymptotic behaviour or stabilization analysis for more general equations, linear or nonlinear. It  all depends a lot on the particular features, and only one or some of the methods work in each case. Other settings of the Heat Equation are briefly discussed in  Section 9 and a longer mention of results for different equations is done in Section 10.
\end{abstract}

\

\

   \medskip

\noindent \it 2010 Mathematics Subject Classification\rm: 35K05, 35K08.


\noindent\it Keywords and phrases\rm: Heat equation, asymptotic behaviour, convergence to the Gaussian solution.

\newpage


\tableofcontents
\normalsize

\newpage

\section{The Cauchy Problem in $\RR^N$ for the Heat Equation}\label{sec.intro}
The classical heat equation (HE) \ $\partial_t u = \Delta u$, is one of the most important objects in the theory of partial differential equations, and it has great relevance in the applied sciences. It has been developed since the seminal work of J. Fourier, 1822, \cite{Fou1822}, to describe phenomena of heat transport and diffusion in many contexts. The great progress of the mathematical theory in these two centuries has had a strong influence not only on PDEs, but also on Probability, Functional Analysis,  as well as  Numerics. It now has well-established connections with other subjects like viscous fluids, differential geometry, image processing or finance. See more on motivation in standard textbooks, like \cite{EvansPDE, SalsaBk} or in the survey paper \cite{VazCIME}.

The heat equation enjoys a well developed theory that has many distinctive features. It can be solved in many settings, suitable for different interests, that lead to quite different results and use different tools. In this paper we will study  the most typical setting, the Cauchy Problem  posed in $\RR^N$, $N\ge 1$:
\begin{equation}\label{eq:CPHE}\tag{CPHE}
\begin{cases}
\begin{aligned}
\partial_t u = \Delta u \qquad\quad\;\, &\text{in } \RR^N\times(0,\infty) \\
u(x,0) = u_0(x) \quad &\text{in } \RR^N,
\end{aligned}
\end{cases}
\end{equation}
with initial data $u_0 \in L^1(\RR^N)$. The equation appears in the theory of Probability and Stochastic Processes as the PDE description of Brownian motion, and in that application the function $u(\cdot,t)$ denotes the probability density of the process at time $t$, hence it must be nonnegative and the total integral $\int u(x,t)\,dx$, often called the mass, must be one. Such restrictions on sign and integral are not needed in the analytical study, and the extra generality is convenient both for the theory and for a number of other applications.

The first question to be addressed by mathematicians is to find appropriate functional spaces to solve this problem and to prove that it is a well-posed problem in the sense of Hadamard's definition (existence, uniqueness and stability). In tackling such a problem we will be lucky enough to find a very explicit representation of the solution by means of the so-called ``fundamental solution'' procedure.

\medskip

\noindent  {\bf Definition.} We call  {\sl fundamental solution of the HE in $\RR^N$} the solution of the equation with initial data a ``Dirac delta'':
\begin{equation}\label{eq:CPHE1}\tag{FS}
\begin{cases}
\begin{aligned}
\partial_t U = \Delta U \qquad\quad\;\, &\text{en } \RR^N\times(0,\infty) \\
U(x,0) = \delta_0(x) \quad\;\, &\text{en } \RR^N.
\end{aligned}
\end{cases}
\end{equation}
In probabilistic terms, we start from a unit point mass distribution  concentrated at the origin $x_0=0$, and we look for the way it expands with time according to the heat equation.

\medskip

\noindent {\bf Exercise 1.} (i) Prove by the method of Fourier Transform that the fundamental solution is given by the formula
\begin{equation}\label{Gaussiana}\tag{GF}
U(x,t) = (4\pi t)^{-\frac{N}{2}} e^{-\frac{|x|^2}{4t}}.
\end{equation}
\noindent {\sl Sketch.}  The Fourier transform of $U$ satisfies $\partial_t \widetilde U = -|\xi|^2\widetilde  U$ with $\widetilde U(\xi,0)=1$, hence $\widetilde U(\xi,t)=e^{-|\xi|^2 t}$. Then, invert the Fourier transform.

(ii) Prove that $ U(x, t) \to \delta_0 (x) $ when $ t \to 0 $ in the sense of distributions, i.\,e.,
\[
\int_{\RR^N}U(x,t)\phi(x)dx \; \to \; <\delta_0,\phi> := \phi(0), \quad \text{for all } \phi \in C_c^{\infty}(\RR^N), \quad \text{ as \ } t \to 0.
\]
(iii) Prove that for every sequence $ t_n \to 0 $ the functional sequence $ U (x, t_n) $ is a $ C^\infty $  approximation of the identity.

\medskip

\begin{center}
\includegraphics[width=3.5cm]{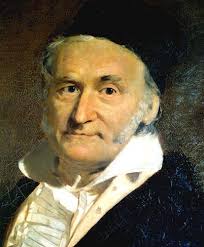}\qquad \qquad
\includegraphics[width=5.0cm]{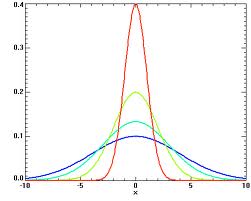}
\end{center}

\centerline{C. F. Gauss and the Gaussian solution}

\medskip

\noindent {\bf Notations and remarks.} This function is known as the Gaussian Function or Gaussian Kernel and we use the notation $ G(x, t) $ or $ G_t (x) $ instead of $U$. We will also use the Gaussian measure $ d\mu = G_tdx $; it is a Radon measure in $ \ren $ with parameter $ t $.  The  Gaussian   kernel is  named after  the famous  German mathematician Carl  Friedrich Gauss, the Prince of Mathematics.

In statistics and probability theory the Gaussian function appears as the density function of the normal distribution. Then   the standard deviation of the distribution is taken as the  parameter in the formula, and we have \ $ 2t = \sigma^2 $. This relation : {\sl time proportional to square space deviation}, is called the Brownian scale and appears in many calculations of evolution processes related to heat propagation or Brownian motion. Finally, the same function appears in Statistical Mechanics with the name of Maxwellian distribution (after James Clerk Maxwell, 1860). See other historical comments at the end of the paper.

A main result in those disciplines is the fact that $G$ appears as the  limiting probability distribution of many discrete processes, according to the famous Central Limit Theorem. The limiting behaviour of the solutions of the Heat Equation is our main concern in this paper, and the Gaussian kernel will play a main role.

\medskip

\noindent  {\bf Exercise 2.} Calculate the Gaussian function by the method of self-similar solutions, $U(x,t)=t^{-\alpha}F(x\,t^{-\beta})$,  plus the condition of mass conservation.

{\sl Hint.} Substitute the self-similar form into the PDE and check that time disappears explicitly from the resulting equation for $F=F(\xi)$ (with $\xi=x\,t^{-\beta}$) if we put $\beta=1/2$. Check that the choice $\alpha=N/2$ means conservation of mass. This selects both similarity exponents. Then write the resulting elliptic equation for $F(\xi)$ assuming that it is a function of $r=|\xi|$ (it is a radial function) to get the ODE
$$
(r^{N-1}F'(r))'+ \frac{1}{2}(r^N F)'=0.
$$
Calculate that \ $\log F(r)= a- r^2/4$, with $a\in\re$, hence $F=Ce^{-|x|^2/4}.$

\medskip

\noindent {\bf Exercise 3.} Write the formula of  Gaussian with  mass $ M $ as
$$
U(x,t;M)=\frac{M}{(4\pi t)^{N/2}}e^{-x^2/4t}
$$
and check that it is dimensionally correct. Dimensional analysis is a powerful tester for correct formulas (if applied well).

\medskip

\noindent {\bf Exercise 4. Gaussian Representation Formula.} Show that if $ u_0 \in C(\RR^N) $ and is bounded, then the convolution
\begin{equation}\label{eq:conv.form}\tag{GRF}
u(x,t) = (U \ast u_0)(x,t) := (4\pi t)^{-\frac{N}{2}} \int_{\RR^N}u_0(y) e^{-\frac{|x-y|^2}{4t}}dy
\end{equation}
is a solution of the heat equation with initial data $u_0$. Check that it is a $C^\infty$ function of $x$ and $t$ for every $x\in\ren$ and $t>0$. Prove that $ u(x,t) \to u_0 (x) $ pointwise when $ t \to 0 $.

(ii) Prove that the initial data is taken uniformly if $ u_0 \in BUC (\RR^N) $, the space of bounded and uniformly continuous functions in $\ren$.

(iii*) Consider whether the boundedness condition of point (i) can be relaxed or eliminated.

\medskip

\noindent {\bf Exercise 5.} Prove that formula \eqref{eq:conv.form} is still valid when the initial data belong to  $L^p(\RR^N)$, $1\le p\le\infty $.

(ii) What is the regularity of the solution?

(iii) In which sense  are the initial data taken?

\medskip

\noindent {\bf Exercise 6.} (i)  Prove that  (GRF) generates a continuous contraction semigroup in the following Banach spaces
$$
X=BUC(\ren), \quad L^2(\ren), \quad L^p(\ren), \quad  1<p<\infty\,,
$$
by means of the definition $S_t:X\to X$, $S_t u_0=u(\cdot,t)$, where $u$ is the solution of the heat equation with initial data $u_0$.

{\sl Sketch.} What you have to check is that, as a function of time $t$ \  $f(t):= S_tu_0$ belongs to $C([0,\infty: X)$, that $S_t\ast S_s=S_{t+s}$, and also $S_tu_0\to u_0$ as $t\to 0$ in $X$. Moreover, contraction means that
$\|S_tu_0\|_X\le \|u_0\|_X$ for every $t>0$ and every $u_0\in X$.

(ii) Study in which sense the semigroups agree (explain if they are the same or not). You may define the {\sl core} of the semigroup by restricting it to a nice and dense set of functions, and then explain how the concept of density is used to recover the complete semigroup in each of the above spaces.

\medskip

\noindent {\sc Notation.} We will often write $u(t)$ instead of $u(\cdot,t)$ for brevity in the hope that no confusion will arise. Hence, $u(t)$ is a function of $x$  for every fixed $t$.

\medskip

\noindent {\bf Exercise 7.} (i)  Prove that for  $p=1$ we have the stronger property of {\bf mass conservation}
\begin{equation}\tag{MC}
\int_{\ren} u(x,t)\,dx=\int_{\ren}u_0(x)\,dx\,,
\end{equation}
for all data in $L^1(\ren)$.

(ii) Note that this law also holds for signed solutions. Show that the version with absolute value is not valid for signed solutions by means of an example (superposition of two Gaussians)

(iii) Prove on the fundamental solution that the  $L^p$-integrals with $p>1$ are not  conserved and check the rate of decay they exhibit. {\sl Terminology:} The $p$ energy is $\int |u|^p \,dx. $ Those $p$-energies are not conserved because they ``undergo dissipation''.

(iv*) The curious reader may want to learn more about that dissipation and how to measure it accurately. Any ideas?

\medskip

\noindent {\bf Exercise 8. Conservation of moments.} (i) We define the first (signed) moment as the vector quantity
\begin{equation}
{\mathcal N}_1(u(t))=\int_{\ren} x\,u(x,t)\,dx\,,
\end{equation}
which in principle evolves with time. Prove that this quantity is  actually conserved in time  for any solution of the HE with initial data such that
$\int (1+|x|)|\,u_0(x)|\,dx$ is finite.

(ii) For data such that $\int |x|^2\,|u_0(x)|\,dx$ is finite, we define the second moment as
the scalar time-dependent quantity
\begin{equation}
{\mathcal N}_2(u(t))=\int_{\ren} |x|^2\,u(x,t)\,dx\,.
\end{equation}
Prove that this quantity is not conserved in time  for solutions of the HE, and in fact
\begin{equation}
{\mathcal N}_2(u(t))={\mathcal N}_2(u_0) +2Nt\,.
\end{equation}
(Recall here that $N$ is the space dimension). These results will be important to understand the finer versions of the asymptotic convergence to the Gaussian. Deeper analysis proves that the only integral quantities conserved by the solutions of the Heat Equation are the mass and the first moment.

\medskip

\noindent {\bf Exercise 9.} (i) Prove the so-called {\bf ultra-contractive estimates}, also called {\sl smoothing effects}. They allow to pass from data in $L^p$ to solutions $u(t)$ in $ L^\infty$ for every $t>0$. More precisely, prove that there exists a constant $C=C(p,N)$ such that for every $u_0 \in L^p(\ren)$
$$
\|u(t)\|_\infty \le C\frac{\|u_0\|_p}{t^{N/2p}}.
$$
We may say that the constant is universal since it does not depend of the particular solution we take.

(ii) Calculate the exponent $N/2p$ using only dimensional calculus.

\medskip

\noindent {\bf Exercise 10. Regularity}. (i) Prove that there exists a universal constant $C=C(p,N)>0$ such that for every $p\ge 1$ and every $t>0$
$$
\|\frac{du(t)}{dt}\|_p\le C\,\frac{\|u_0\|_p}{t}\,.
$$

(ii) Prove that for nonnegative solutions we have a stronger pointwise inequality
$$
\frac{\partial u}{\partial t} \ge -C_1  \frac{u}{t}\quad \mbox{with} \quad C_1=N/2\,.
$$
Prove that this constant is optimal. {\sl Hint.} Check the properties first for the fundamental solution, then use convolution.

(iii)  Prove that for every  $ t > 0 $ the $x$-derivatives of $u$ satisfy
\begin{equation*} \label{ESTIMACIONDERIVADAS}
 \left|\frac{\partial}{\partial x_i} u(x,t)\right| \leq C \frac{|| u_0||_{L^{\infty}}}{t^{\frac{1}{2}}}.
\end{equation*}
Obtain a similar estimate for data in  $ L^p\left(\RR^N \right) $.

(iv) Check that the dimensions are all correct in these formulas.

\medskip

\noindent {\bf Exercise 11. Gaussian Representation Formula with measure data.} (i) Prove that the representation formula can be used with any bounded Radon measure as initial data, in the form
\begin{equation}\label{eq:conv.form.m}\tag{GRFm}
u(x,t) = (G_t \ast \mu)(x,t) := (4\pi t)^{-\frac{N}{2}} \int_{\ren} e^{-\frac{|x-y|^2}{4t}}d\mu(y)\,.
\end{equation}

(ii) Prove that it produces a classical solution of the HE for all $t>0$ and that the a priori estimates hold with the $L^1$ norm of $u_0$ replaced by the total mass of the measure,
$|\mu(\ren)|=\int |d\mu|$.

(iii) Prove that the initial data are taken in the weak sense of measures.

(iv) Observe that the fundamental solution $U=G_t$ falls into this class.

\medskip

\noindent {\bf Exercise 12.} (i) Find explicit solutions of polynomial type of heat equation. {\sl Suggestion:} try the  solution types
$$
u = A_0 (x), \quad u = A_0 (x) + A_1 (x) t, \quad u = A_0 (x) + A_1 (x) t + A_2 (x) t ^ 2, \dots
$$
Find all the solutions of the first two types.
Maybe the best known solutions in this class are $u=1$ and $u=|x|^2+2Nt$.

(ii) Prove that they also satisfy the representation formula, even though they are not integrable or bounded.

\medskip

The general theory says that the (GRF) holds for all initial data which are locally bounded measures $\mu$ with a weighted integrability condition
$$
\int_{\ren} e^{-c |x|^2}|d\mu(x)|<\infty \quad \mbox{for some} \ c>0\,,
$$
see \cite{Widder}, but we  will not use that sharp result in these notes. We will only point out
the following example of blow-up solution with very large growth as $|x|\to\infty$
\begin{equation}\label{bu.sol}
U_b(x,t)=C\,(T-t)^{-N/2}e^{|x|^2/4(T-t)}.
\end{equation}

\noindent {\bf Exercise 13.} Check that $U_b$ is a classical solution of the heat equation for $-\infty<t<T$ and blows up everywhere in $x$ as $t\to T$.
\medskip

\medskip

\noindent {\bf Exercise 14.} Construction of new solutions.

(i) Show that if $u (x, t) $ is a solution of the heat equation, so are \ $ u_1 = \partial_t u $, and $ u_2 = \partial_{x_i} u $.

(ii) Show that if $u (x, t) $ is a solution of the heat equation in 1D then so is also
 $ u_3 = \int_{- \infty}^x u(s, t) \, ds $.

(iii) Prove that if $u (x, t) $ is a solution of the heat equation so is also
$$
v(x, t) = A u(Bx, B^2t)
$$
for any choice of the parameters $ A $ and $ B $. This is called the scaling property.
$ A $ and $ B $ need not be positive.

(v*) Show that if $u (x, t) $ is a solution of the heat equation, then so is also
$$
v(x, t) = x u_x + 2t u_t\,.
$$
This is a 1D notation. But the result is valid in any dimension if the notation is correctly interpreted.

%
\section{Asymptotic convergence to the Gaussian}\label{main.result}

The main result on the asymptotic behaviour of general integrable solutions of the heat equation consists in proving that they look increasingly like the fundamental solution.   Since this solution goes to zero uniformly with time, the estimate of the convergence has to take into account that fact and compensate for it. This happens by considering a renormalized error that divides the standard error in some norm by the size of the Gaussian solution $U(t)=G_t$ in the same norm. For instance, in the case of the sup norm we know that
$$
\|G_t\|_{L^\infty(\ren)}=Ct^{-N/2},\qquad \|G_t\|_{L^1(\ren)}=1.
$$
This is the basic result we want to prove

\begin{theorem} \label{main.convthm} Let $u_0\in L^1(\ren)$ and let  $\int u_0(x)dx=M$ be its mass. Then the solution $u(t)=u(\cdot,t)$  of the HE in the whole space ends up by looking like $M$ times the fundamental solution  $U(t)=G_t$ in the sense that
\begin{equation}\label{cr.l1}
\lim_{t\to\infty} \|u(t)-M G_t\|_1\to 0
\end{equation}
and also that
\begin{equation}\label{cr.linf}
\lim_{t\to\infty} t^{N/2}\|u(t)-M G_t\|_\infty\to 0\,.
\end{equation}
By interpolation we get the convergence result for all $L^p$ norms
\begin{equation}\label{cr.lp}
\lim_{t\to\infty} t^{N(p-1)/2p}\|u(t)-M G_t\|_{L^p(\ren)}\to 0\,.
\end{equation}
for all $1\le p\le \infty$. 
\end{theorem}

 We add important information to this result in a series of remarks.

\noindent $\bullet$  First, one comment about the spatial domain. The fact that we are working in the whole space is crucial for the result of Theorem \ref{main.convthm}. The behaviour of the solutions of the heat equation posed in a bounded domain with different kinds of boundary conditions is also known, and the asymptotic behaviour {\sl does not follow the Gaussian pattern.}

\noindent $\bullet$  As we will see below, convergence to the Gaussian happens on the condition that the data belong to the class of integrable functions, that can be extended without problem to bounded Radon measures. This is actually no news since the theory says that the solution corresponding to an initial measure $\mu\in {\mathcal M}(\ren)$ is integrable and bounded for any positive time, so we may change the origin of time and make the assumption of integrable and bounded data. But we point out the some of the proofs work directly for measures without any problem.

\noindent $\bullet$ We recall that it is usually assumed that $u\ge 0$ on physical grounds but such assumption is not at all needed for the analytical study of this paper. Thus, the basic result holds also for signed solutions even if the total integral is negative, $M\le 0$. There is no change in the proofs. We may also put $M=1$ by linearity as long as $M\ne 0$.

\noindent $\bullet$   $M=0$ is a special case that deserves attention: even if Theorem \eqref{main.convthm} is true, the statement does not imply that the solution looks asymptotically like a Gaussian; to be more precise, it only says that the previous first-order approximation disappears. If we want more precise details about what the solution looks like, we have to search further to identify the terms that may give us the size and shape of such a solution. This question will be addressed below. For the  moment let us point out that differentiation in $x_i$ of the Heat Kernel produces a new solution with zero integral
$$
u_i(x,t)=\partial_{x_i} G_t(x)= C\,t^{-(N+2)/2}x_i\,e^{-|x|^2/4t}\,,
$$
to which we can apply the above comments. In particular,
$$
t^{N/2}u_i(x,t)= O(t^{-1/2})\,,
$$
where $O(\cdot)$ is the Landau $O$-notation for orders of magnitude.

\noindent $\bullet$ It must be stressed that convergence to the Gaussian {\sl  does not hold for other data}. Maybe  the simplest example of solution that does not approach the Gaussian is given by any non zero constant solution, but it could  be objected that $L^\infty(\ren)$ is very far from $L^1(\ren)$. Actually, the same happens for all $L^p(\ren)$ spaces, $p>1$. Indeed, a simple argument based on approximation and comparison shows that for any $u_0\ge 0$ with $\int u_0(x)\,dx=+\infty$ we have
$$
\lim_{t\to\infty} t^{N/2}u(x,t)=+\infty
$$
everywhere in $x\in \ren$ (and the divergence is locally uniform).

\noindent $\bullet$  The way different classes of non-integrable solutions actually behave for large time is an interesting question that we will not address here. Thus, the reader may prove using the convolution formula that
for locally integrable data that converge to a constant $C$ as $|x|\to\infty$ the solution $u(x,t)$ stabilizes  to that constant as $t\to\infty$. Taking growing data may produce solutions that tend to infinity with time, like the 1D family of travelling waves
$$
U_{TW}(x,t)=C\,e^{c^2t+cx}
$$
defined for real  constants $C,c\ne 0$. A more extreme case is the blow-up solution $U_b(x,t)$ of formula \eqref{bu.sol} that not only does not stay bounded with time, it even blows up in finite time.

\noindent $\bullet$ About the three convergence results of the Theorem, it is clear that \eqref{cr.lp} follows from \eqref{cr.l1} and \eqref{cr.linf}. Now, what is interesting is that \eqref{cr.linf} follows from \eqref{cr.l1} and the smoothing effect (Exercise 9). We ask the reader to prove this fact. Hint: use \eqref{cr.l1} between $0$ and $t/2$ and then the  smoothing effect between $t/2$ and $t$.

\noindent $\bullet$  It is interesting to note that the convergence in $L^1$ norm  of formula \eqref{cr.l1} can be formulated without mention to the existence of  a fundamental solution in the following  form:

\noindent {\bf Alternative Theorem.} {\sl Let $u(t)$ and $v(t)$ be any two solutions of the Cauchy problem for the HE in the whole space, and let us assume that their initial data satisfy $\int u_0\,dx=\int v_0\,dx$. Then we have}
\begin{equation}\label{cr.l1.mod}
\lim_{t\to\infty} \|u(t)-v(t)\|_1\to 0\,.
\end{equation}
We could use a similar approach for the $L^\infty$ estimate but some Gaussian information appears in the weight $t^{N/2}$. The alternative approach has been first remarked by specialists in stochastic processes of Brownian type, and it is known as a {\sl mixing property}. It has been generalized to many variants of the heat equation.

%
\section{Proof of convergence by scaling}\label{sec.scaling}

There are many approaches to the proof of the main result, and we will show some of the best known  below. They are interesting for their possible extension to similar asymptotic results for other equations, both linear and nonlinear, see Sections \ref{sec.appl1} and \ref{sec.appl2}. The first proof we give is based on scaling arguments. In the proposed method the proof is divided into 5 steps.

\medskip

\noindent {\sc Step 1. \sl Scaling transformation.} It is easy to check that the HE is invariant under the following one-parameter family of transformations ${\mathcal T}_k $ defined on space-time functions by the formula
$$
{\mathcal T}_k u(x,t)= k^{N} u (kx,k^2t), \qquad k>0,
$$
see also Exercise 14 (iii). This transformation maps solutions $u$ into new solutions $u_k={\mathcal T}_k  u$, called rescalings of the original solution. It also conserves the mass of the solutions
$$
\int_{\ren} {\mathcal T}_k u(x,t)\,dx=\int_{\ren} u(y,k^2t)\,dy=\mbox{\rm constant}.
$$

\medskip

\noindent {\sc Step 2. \sl Uniform estimates.} The estimates proved on the general solutions in Section \ref{sec.intro} imply that the whole family $u_k$ is uniformly  bounded in space and time if time is not so small: $x\in\ren$ and $t\ge t_1>0$ (cf. Exercise 9). They also have uniform estimates on all derivatives under the same restriction on time (cf. Exercise 10).

\medskip

\noindent {\sc Step 3. \sl Limit problem.} We can now use functional compactness and pass to the limit $k\to \infty$ along suitable subsequences to obtain a function $\hat U(x,t)$ that satisfies the same estimates mentioned above and is a weak solution of the HE in $Q=\ren\times (0,\infty)$.  By the estimates on derivatives the solution is classical. The convergence takes place locally in $Q$ in the sup norm for the functions and their space and time derivatives.

\medskip

\noindent {\sc Step 4. \sl Identifying the limit.} We now have to check that the limit solution $\hat U(x,t)$ has the same mass $M$ as the sequence $u_k$ and that is takes the Dirac delta as initial data. By uniqueness for solutions with measure data (which we accept as part of the theory) we will then conclude that $\hat U(x,t)$ is just a fundamental solution, $\hat U(x,t)=M\,G_t$.

(i) The proof of this step is best done when $u_0$ is  nonnegative, compactly supported and bounded, since in that case the solution is bounded above by constant times the fundamental solution, $|u(x,t)|\le C\,G(x,t+1)$. Applying the transformation we get
$$
|u_k(x,t)|\le C\,{\mathcal T}_k (G(x,t+1))= C\, G(x,t+ k^{-2}).
$$
This means a uniform control from above of the mass of all the tail mass of all the solutions $u_k$, by which we mean the mass lying in exterior regions of space.  Such control allows to avoid the loss of mass that could occur in the limit by Fatou's Theorem. Hence, the mass of $\hat U$ is the same, i.e., $M$.

The convergence of $\hat U$ to $M\,\delta(x)$ as $t\to 0$ happens because $u_k(x,0)\to M\,\delta(x)$, and the previous tail analysis shows that $\hat U$ takes zero initial values for $x\ne 0$.

(ii) To recover the same result for a general $u_0 \ge $ we use approximation by data as above and then the $L^1$ contraction of the heat semigroup. For signed solutions, separate the positive and negative parts of the data and solve separately.

\medskip

\noindent {\sc Step 5. \sl Recovering the result.} (i) We now use the convergence of the rescaled solutions at a fixed time, say $t=1$,
$$
\lim_{k\to\infty}|u_k(x,1)-MG_1(x)|\to 0.
$$
This convergence takes place locally in $\ren$ in all $L^p$ norms. By the tail analysis, it also happens in $L^1(\ren)$. But since the derivatives are also uniformly bounded we have also an $L^\infty$ estimate.

\medskip

(ii) We now use the meaning to the transformation ${\mathcal T}_k$ and write $k=t_1^{1/2}$ to get
$$
\lim_{t_1\to\infty} |t_1^{N/2}u(xt_1^{1/2},t_1)-MG_1(x)|\to 0.
$$
But this is just the result we wanted to prove after writing $x=y\,t_1^{-1/2}$ and observing that
$G_{t_1}(y)=t_1^{-N/2}G_1(x)$. Similarly for the $L^1$ norm.

\medskip

This 5-step proof is taken from paper  \cite{KVplap}, where it was applied to the $p$-Laplacian equation. See whole details for the porous medium equation in the book \cite{VazBook}. It has had further applicability.

\medskip

\noindent {\bf Exercise 15.} Fill in the details of the above proof.

%
\section{Asymptotic convergence to the Gaussian via representation}\label{sec.conv1}

The second proof we give of the main result, Theorem \ref{main.convthm}, is based on the examination of the error in terms of the representation formula. This is a very direct approach, but it needs a previous step,  whereby the proof is done under the further restriction that the data have a finite first moment. Then the convergence result is more precise and quantitative. This particular case has an interest in itself since it shows the importance of controlling the first moment of a mass distribution. We already know that the first moment is associated to a conserved quantity (see Exercise 8).

\begin{theorem} \label{thm.conv.2} Under the assumptions that $u_0\in L^1(\ren)$ and that the first absolute moment is finite
\begin{equation}
{\mathcal N}_1=\int |u_0(y)y|\,dy <    \infty\,,
\end{equation}
we  get the convergence
\begin{equation}\label{conv.fm}
\displaystyle t^{N/2}|u(x,t)-MG_t(x)|\le C{\mathcal N}_1t^{-1/2}\,,
\end{equation}
as well as
\begin{equation}
\displaystyle \|u(x,t)-MG_t(x)\|_{L^1(\ren)}\le  C{\mathcal N}_1t^{-1/2}\,.
\end{equation}
The rate $O(t^{-1/2})$ is optimal under such assumptions.
\end{theorem}

\noindent {\sl Proof.} (i) We may perform the proof under the further restriction that $u_0\ge 0$. For a signed solution we must only separate the positive and negative parts of the data and apply the results to both partial solutions.

(ii) Let us do first the sup convergence. We  have
$$
\begin{array}{c}
\displaystyle u(x,t)-MG_t(x)= \int u_0(y)G_t(x-y)\,dy- G_t(x)\int u_0(y)\,dy\\[4pt]
\displaystyle =\int u_0(y)(G_t(x-y)-G_t(x))\,dy\\
\displaystyle =\int u_0(y)\left(\int_0^1\partial_s (G_t(x-sy))\,ds\right)\,dy\\[4pt]
\displaystyle =  C t^{-N/2}
\int dy \int_0^1 ds \,u_0(y) \langle y, \frac{x-sy}{2t}\rangle e^{-|x-sy|^2/4t}\\[4pt]
\end{array}
$$
with $C=(4\pi)^{-N/2}$.
Consider the piece of the integrand of the form
$$
f=\frac{x-sy}{t^{1/2}}e^{-|x-sy|^2/4t}= \xi e^{-|\xi|^2/4}\,,
$$
where we have put $\xi=(x-sy)/t^{1/2}$. We observe that the vector function $f$ is bounded
by a numerical constant, hence
$$
\displaystyle |u(x,t)-MG_t(x)|\le  C_1t^{-(N+1)/2}\int |u_0(y)y|\,dy\,.
$$
 Taking into account that $ G_t$ is of order $t^{-N/2}$ in sup norm, we write the result as \eqref{conv.fm}, and  $C=C(N)$ is a universal constant.

 \medskip

 (iii) For the $L^1$ convergence we start in the same way and arrive at
$$
\displaystyle u(x,t)-MG_t(x)=  C \,t^{-N/2}
\int dy \int_0^1 ds \,u_0(y) \langle y, \frac{x-sy}{2t}\rangle e^{-(x-sy)^2/4t}\,.
$$
Now we integrate to get
$$
\begin{array}{c}
\displaystyle \|u(x,t)-MG_t(x)\|_1\le C \, t^{-N/2}
\int_{\ren} dx\int_{\ren} dy \int_0^1 ds \,|u_0(y) y| \frac{|x-sy|}{2t}e^{-|x-sy|^2/4t}\\[6pt]
\displaystyle = C \,\int_{\ren} dy \int_0^1 ds \,|u_0(y) y|t^{-1/2} (\int_{\ren } t^{-N/2}\frac{|x-sy|}{t^{1/2}}e^{-|x-sy|^2/4t} \,dx)
\end{array}
$$
With the change of variables $x-sy=t^{1/2}\xi$ we already know that the last integral is a constant independent of $u$, hence the formula for the $L^1$ error.
Note that now we are speaking of masses and we do not need any renormalization time factor.

(iv) We ask the reader to prove the optimality as an exercise.
\qed

\medskip

\noindent {\bf Exercise 16.} Take as solution the Gaussian after a space displacement, $u(x,t)= G(x+h,t)$, and find the convergence rate to be exactly $O(t^{-1/2})$. This is just a calculus exercise but attention to details is needed. {\sl Hint}: Write
$$
G(x+h,t)=G_t(x)+ h\partial_x G_t(x)+ \frac{h^2}2 D^2_xG_t(\xi)
$$
(where $\xi=x+sh$, $0<s<1$), and check that
$$
\partial_x G_t(x)=Ct^{-(N+1)/2}\,{\xi}e^{-\xi^2/4}=O(t^{-(N+1)/2})\,,
$$
uniformly in $x$,  and $D^2_xG_t(\xi)= O(t^{-(N+2)/2})$ uniformly in $x$. This exercise shows that the term $h\partial_x G_t(x)$ is the {\sl precise corrector} with relative error $O(t^{-1/2})$. We could continue the analysis by expanding in Taylor series with further terms, see below.   Exact correctors and longer expansions for general solutions will be done later in Section \ref{sec.spec} by the methods of Functional Analysis.

\medskip

\noindent {\bf Remark.} The factor $t^{N/2}$ is the appropriate weight to consider relative error.
We point out that the more precise relative error formula
\begin{equation}
\epsilon_{rel}(t)=\frac{|u(x,t)-G_t(x)|}{G_t(x)}
\end{equation}
does not admit a sup bound, as can be observed by choosing $u(x,t)=G_t(x-h)$ for some constant $h$ since then
$$
\epsilon_{rel}(t)=|e^{xh/2t}e^{-h^2/4t}-1|\,,
$$
which is not even bounded. It is then quite good that our weaker form does admit a good estimate. Same happens for the $L^1$ norm. This comment wants to show that error calculations with Gaussians are delicate because of the tail (i.e., the behaviour for large $|x|$.

\medskip

\noindent $\bullet$ {\sl Proof of Theorem \ref{main.convthm} in this approach}. Given an initial function $u_0\in L^1(\ren)$ without any assumption on the first moment, we argue by approximation plus the triangular inequality. In the end we get a convergence result, but less precise. To quote, if  $u_0$ is integrable with integral $M$ and let us fix an error  $\delta>0$. First, we find an approximation  $u_{01}$ with compact support and such that
$$
\|u_0-u_{01}\|_1<\delta.
$$
Due to the already mentioned effect $L^1\to L^\infty$, and applied to the solution $u-u_1$, we know that for all  $t>0$
$$
\|u(t)-u_{1}(t)\|_\infty<C\,\delta\, t^{-N/2}\,.
$$
On the other hand, we have just proved that for data with finite moment:
$$
\|u_1(t)-M_1G_t\|_\infty\le C {\mathcal N}_{1\delta} \,t^{-(N+1)/2}\,.
$$
In this way, for sufficiently large $t$  (depending on $\delta$) we have
$$
 t^{N/2}\|u_1(t)-M_1G_t\|_\infty\le C \delta
$$
with $ C$ a universal constant. Next se recall that  $|M-M_1|\le \delta$ as well as $G_t\le Ct^{-N/2}$. Using the triangular inequality we arrive at
\begin{equation}
\lim_{t\to\infty}t^{N/2}\|u(t)-MG_t\|_\infty=0\,,
\end{equation}
which ends the proof. \qed

\noindent {\bf Remark}. In the general conditions of Theorem \ref{main.convthm} we still obtain convergence, but we no longer obtain a rate. In fact, we show next that no convergence speed can be found without further information on the data other than integrability of  $ u_0 $.

\medskip


 \hiddensubsection{No explicit rates for general data. Counterexample}

Let us explain how the lack of a rate for the whole class of $L^1$ functions is shown. Given any decreasing and positive rate function $\phi$ such that $\phi(t)\to 0$ as $t\to\infty$ we construct a modification of the Gaussian kernel that produces a solution with the same mass $M=1$ and such that it satisfies a lower bound for the error of the form
$$
 t^{N/2}\|u(x,t)- G_t(x)\|_\infty\ge n \phi(t_n)
$$
at a sequence of times $t_n\to\infty$ to be chosen.

\noindent {\bf Construcion}. The idea is to find  a choice of small masses $m_1,m_2,\dots$ with $\sum_n m_n=\delta< 1$, and locations $x_n$ with $|x_n|=r_n\to \infty$ and consider the solution
$$
u(x,t)=(1-\delta)G_t(x)+ \sum_{n=1}^\infty m_nG_t(x-x_n)\,.
$$
Let us be precise.  The error $u(x,t)-G_t(x)$  is calculated at $x=0$ as
$$
(4\pi t)^{N/2}|u(0,t)-G_t(0)|=|\delta-\sum_{n=1}^\infty m_n e^{-x_n^2/4t}|=
\sum_{n=1}^\infty m_n (1-e^{-x_n^2/4t})\,.
$$
 Put $m_n=2^{-n}$ (any other summable series will do). Choose iteratively $t_n$ and $x_n$ as follows. Given choices for the steps $1,2,\dots, n-1$, pick $t_n$ to be much larger than $t_{n-1}$ and such that $\phi(t_n)\le m_n/2n$. This is where we use the fact that $\phi(t)$ tends to zero, even if it may decrease in a very slow way. Choose now $|x_n|=r_n$ so large that \ $e^{-r_n^2/4t_n}< 1/2$. Essentially, the mass has to be displaced at distance equal or larger than $O(t_n^{1/2})$. Then,
$$
m_n(1-e^{-x_n^2/4t_n})\ge n\phi(t_n)\,.
$$

Let us make some further practical calculations (with no precise scope in mind): Let for example $\phi(t)=t^{-\epsilon}$. Choose $m_n=2^{-n}$. Then
$t_n ^{\epsilon}\ge 2^{n}$,
$$
t_n\sim 2^{n /\epsilon}, \qquad r_n\sim 2^{n/2\epsilon}
$$
and the mass in the outer region $\{|x|\ge r_n\}  $  is approximately $2^{-n}$; hence
such a mass is $M(r)\sim r^{-2\epsilon}\,,$ which is not so small if $\epsilon\to 0$.
\normalcolor

 \hiddensubsection{Infinite propagation in space }

 The representation formula immediately shows that a solution $u$ corresponding to nonnegative initial data will be
 strictly positive at all points $x\in\ren$ for any time $t>0$. The infinite speed of propagation of the heat equation with the instantaneous formation  of a thin tail at infinity is considered an un-physical property  by many authors,  one of the not many drawbacks of this wonderful equation. However, it is essential to the equation and creates some curious effects.

 \medskip

 \noindent $\bullet$ {\large \bf Spatial tails for positive solutions.} Let us examine the precise form of the tail in the simplest case where $u_0$ has compact support. For simplicity and w.l.o.g. we assume that $u_0$ is supported in the ball $B_R$ of radius $R>0$ centered at $x=0$. Take $x_1\in \ren$ such that $|x_1|=r_1>1$, say $x_1=r_1 {\bf e}_1$. Then the representation formula implies that a bound from above  is obtained by displacing all the mass to the nearest point to $x_1$ inside $\overline{B_R}$, which is $x_0'=R {\bf e}_1$, and we get
  the upper bound
 $$
 u(x_1,t)\le M\,G_t(x_1-x'_0)=\frac{M}{(4\pi t)^{N/2}}e^{-(|x_1|-R)^2/ 4t}   \,.
 $$
In the opposite direction, moving all the mass to $x_0''=-R {\bf e}_1$ we get the lower bound
$$
u(x_1,t)\ge M\,G_t(x_1-x_0'')=\frac{M}{(4\pi t)^{N/2}}e^{-(|x_1|+R)^2/ 4t}   \,.
 $$
 \normalcolor
 Both estimates are clearly optimal in this context.
Since the equation is invariant under rotations the result holds for all $x$ such that $|x|>R$ instead of our restricted choice of $x_1$. Moreover, we see that the ratio of both estimates tends to infinity as $|x|\to\infty$. It is therefore convenient to take logarithms, and then we easily get a general formula where $B_R$ is centered at $x_c$ with $ |x_c|=R_c$.

\begin{proposition} \label{prop.tails} Let $u_0$ be supported in the ball $B_R(x_c)$ and let $M=\int u_0(x)\,dx$. Then for every $t>0$ and every  $x\not\in B_R(x_c)$ we have
\begin{equation}
\frac{N}2 \log(4\pi t)+\frac1{4t}(|x-x_c|-R)^2 \le -\log (u(x,t)/M)
\le \frac{N}2 \log(4\pi t)+\frac1{4t}(|x-x_c|+R)^2   \,.
\end{equation}
 It follows that
 \begin{equation}
 \lim_{|x|\to\infty}\frac{\log (u(x,t)/M)}{|x|^2}= -\frac1{4t}.
 \end{equation}
 \end{proposition}

 We conclude that in first approximation the tail at infinity of all solutions with compactly supported initial data is universal in shape and depends only on the mass of the data and time. Of course, the second term in the expansion depends also on the radii $R_c$ and $R$.

 Another observation is that the asymptotic space behaviour allows to calculate the time elapsed since the solution had compact support (if we already know the mass $M$).

 Let us also remark that solutions with more general nonnegative data can have other type of tails and we invite the reader to calculate some of them, both for integrable and non-integrable data. Here is an example that decays like a simple exponential.

\noindent {\bf Exercise 17.} Consider the heat equation in 1D. Show that when $u_0$ is integrable, positive and bounded and $u_0(x)=e^{-x}$ for $x\ge 0$, then for all times
$$
\lim_{x\to\infty} e^x u(x,t)=  e^t.
$$
{\sl Sketch:} Putting $u_0(x)=0$ for $x<0$ for simplicity, use the representation formula to write
$$
e^{x-t}u(x,t)=\frac1{(4\pi t)^{1/2}}\int_{-\infty}^x e^{-(y-2t)^2/4t}dy
$$
Let then $x\to \infty$. You may also use the explicit solution $U(x,t)=e^{t-x+c}$.

 In any case the behaviour described in Proposition \ref{prop.tails}  is the \sl minimal \rm one for nonnegative solutions.

 \medskip


 \noindent $\bullet$ {\large \bf Signed solutions do not become everywhere positive.}\\
The square exponential tail behaviour of the fundamental solution has more curious consequences. Thus, if a signed solution $u(x,t)$ has initial data such that the mass of the positive part $M_+=\int u_0^+(x)dx$  is larger than the mass of the negative part $M_-=\int u_0^-(x)dx$,  then we know that it  converges as $t\to\infty$ to the positive
Gaussian $M G_t$ with $M= M_+-M_- > 0 $ so that
\begin{equation}
 \lim_{t\to\infty} t^{N/2} u(x,t)=(4\pi)^{-N/2}\,M>0.
\end{equation}
for every $x\in\ren$.  It would be natural to expect that when $M_+$ is much larger than $M_-$ then the solution $u$ is indeed positive everywhere for large enough, finite times. Now, this is true in most of the space because of the previous convergence, but it is not true in all the space.

\medskip

\noindent {\bf Exercise 18.} (i) Take $N=1$. Show that for the choice $u_0(x)=\delta(x)-\ve \,\delta(x-1)$ with $0<\ve<1$ (a combination of delta functions) the solution $u$ is positive on the left of a line $x=r(t)$ and negative for $x>r(t)$. Show that for large times $r(t)\sim 2\log(1/\ve) t$.

{\sl Remark.} We see that  the problem arises at the far away tail. Of course, $u(x,t)$ becomes positive  for large times at all points located at or less than the typical distance, $|x|\le C t^{1/2}$.

(ii) Show that a similar result is true for integrable data if $u_0$ is positive for $x<0$ and negative for $x>1$, and zero in the middle. Show in particular that $u(x,t)<0$ if $x> 2t\log(1/\ve) +1/2$.

(iii)  State  similar results in several dimensions.


\section{Improved convergence for distributions with second moment}
\label{sec.conv2}

Better convergence rates can be obtained by asking a better decay at infinity of $u_0$. The technical condition we use is having a finite second moment, a condition that is very popular in the literature. In probability this is known as having a finite variation. The motivating example is described next.

\medskip

\noindent {\bf Exercise 19.}  Consider a  time displacement of the fundamental solution and show that $u(x,t)= G(x,t+t_0)$ has the precise convergence rate $O(1/t)$ towards $G(x,t)$.

The result we prove is as follows.

\begin{theorem}\label{thm.conv.3} Under the assumptions that $u_0\in L^1(\ren)$ and that the signed first moment
\begin{equation}
{\mathcal N}_{1,i}=\int_{\ren} y_i\,u_0(y)\,dy
\end{equation}
is finite (for all coordinates), as well as  the second moment:
\begin{equation}
{\mathcal N}_2=\int_{\ren} y^2\,|u_0|(y)\,dy<\infty,
\end{equation}
we  get the convergence
\begin{equation}\label{conv.fm2}
\displaystyle t^{N/2}|u(x,t)-M\,G_t(x)+\sum_i {\mathcal N}_{1,i}\,\partial_{x_i}G_t(x)|\le
C{\mathcal N}_2\,t^{-1}\,,
\end{equation}
and
\begin{equation}\label{imprvdcon1}
\displaystyle \|u(x,t)-M\,G_t(x)+\sum_i {\mathcal N}_{1,i}\,\partial_{x_i}G_t(x)\|_{L^1(\ren)}\le C\,{\mathcal N}_2 \,t^{-1}\,.
\end{equation}
The rate $O(t^{-1})$ is optimal under such assumptions.
\end{theorem}

\noindent {\bf Remark.} The signed first moment is called in Mechanics the center of mass, for $M\ne 1$, $M\ne 0$, we use the formula
\begin{equation}\label{com}
x_c=\frac1{M}\int x\, u_0(x)dx\,.
\end{equation}
In probability it is  the average location of the sample and $M=1$. When $x_c$ is finite and $M\ne 0$, the center of mass can be reduced to zero by just a displacement of the spatial axis. This very much simplifies formulas \eqref{conv.fm2}, \eqref{imprvdcon1}, see below.

    \medskip

\noindent {\sl Proof.} (i) Starting as in Theorem \ref{thm.conv.2} we arrive at the formula
$$
D:=\displaystyle u(x,t)-MG_t(x) = (4\pi t)^{-N/2} \int_{\RR^N} \left(e^{-|x-y|^2/4t}-e^{-|x|^2/4t}\right)\,u_0(y)\,dy.
$$
Let us consider the 1D function  $f(s)=e^{-|x-sy|^2/4t}$, $s\in\re$, and let us use the Taylor formula
$$
f(1)=f(0)+ f'(0)+\int_0^1 f''(s)(1-s)\,ds\,.
$$
We have
$$
f'(s)=\frac1{2t}\langle y, x-sy\rangle e^{-|x-sy|^2/4t}, \quad
f''(s)=(-\frac1{2t}|y|^2 +\frac1{4t^2}|\langle y, x-sy\rangle|^2)e^{-|x-sy|^2/4t}\,.
$$
Using these results, we get $D=D_1+D_2$, where
$$
D_1= (4\pi t)^{-(N+1)/2} \int u_0(y)\langle y, \frac{x}{2t^{1/2}}\rangle e^{-|x^2|^2/4t}\,dy
= -\sum_i \left(\int y_i u_0(y)\,dy \right) \partial_i G_t(x)\,.
$$
On the other hand,  putting $\xi=(x-sy)/2\sqrt{t}$ and $\widehat y= y/|y|$, we also have
$$
D_2= (4\pi t)^{-(N+2)/2} \int_{\ren} \int_0^1 |y|^2 u_0(y)\left(-\frac12 + (\langle \widehat y, \xi \rangle)^2 \right)  e^{-\xi^2}\,dyds.
$$
Since the factor dependent on $\xi$ is uniformly bounded for all $\xi$ we have
$$
|D_2|\le C(N)\left(\int_{\ren} u_0(y)\, |y|^2\,dy)\right) t^{-(N+2)/2}.
$$
This proves the result.

(ii) We leave to the reader to prove the corresponding statement in $L^1$ norm.

(iii) Optimality of the rate follows from Exercise 19. \qed

\medskip

\noindent {\bf Exercise 20.} Give examples of well-known probability  distributions for which
the first moment is finite or infinite. Same for the second moment. Try with examples of the form
$$
f(x)=C(1+|x|^2)^{-a}\,.
$$
For $a=1$ we get the well-known Cauchy distribution that is integrable only in 1D.\\
{\sl Answers.} Integrable $2a>N$, 1st moment $2a>N+1$, 2nd moment $2a>N+2$.

\medskip

\noindent {\bf Reformulation of the result.} It is well-known that the first moment can be eliminated to moving the origin of coordinates to the center of mass $x_c$ defined in formula \eqref{com} when $M\ne 0$. If we do that and then apply Theorem \ref{thm.conv.3} we get the following result

\begin{corollary} Under the assumptions of Theorem \ref{thm.conv.3}  we have
\begin{equation}
\left\{
\begin{array}{l}
\displaystyle |u(x,t) - M G_t(x-x_c)|\le C\,{\mathcal N}_2^*(u_0)\,t^{-(N+2)/2}), \\[6pt]
\displaystyle \|u(x,t) - M G_t(x-x_c)\|_1\le C\,{\mathcal N}_2^*(u_0)\,t^{-1})\,,
\end{array}
\right.
\end{equation}
where ${\mathcal N}_2^*(u_0)= \int_{\ren} |x-x_c|^2 u_0(x)\,dx$ is the centered second moment.
\end{corollary}

\medskip

\noindent {\bf Higher development result.} A continuation of this method into stricter convergence rates using higher moments can be done by using further terms in the Taylor series development. We will not do it but only quote the statement that can be found in \cite{DuoZZ}.

\noindent {\bf Theorem 4 \cite{DuoZZ}} {\sl Let $G(x, t) $ be the heat kernel. For any $1 \le p \le N/(N - 1) $ and $k \ge 0$ an integer the solution of initial value problem for the heat equation satisfies:
$$
\|u(x,t)-\sum_{|\alpha |\le k-j}\frac{(-1)^{|\alpha|}}{\alpha!}\left(\int f(x)x^\alpha\, dx\right)\,D^\alpha G(x,t)\|_p\le C_k t^{-(k+1)/2}\||x|^{k+1}f(x)\|_p,
$$
for any initial data $f \in  L^1(\ren, 1 + |x|^k)$ such that $ |x|^{k+1}f\in L^p(\ren).$}

Let us note that their proof is different from the previous ones and interesting. As a conclusion, we also have a result about convergence in the first moment norm.

\begin{theorem}\label{thm.conv.4} Under the assumptions that $u_0\in L^1(\ren)$ and that the first, second and third moments are finite, we  get the convergence
\begin{equation}\label{imprvdcon1b}
\displaystyle \|u(x,t)-M\,G_t(x)+\sum_i {\mathcal N}_{1,i}\,\partial_{x_i}G_t(x)\|_{L^1(|x|dx)}\le C \,t^{-1/2}\,.
\end{equation}
The constant $C>0$ depends on $u_0$ through the second and third moments.
\end{theorem}

   \medskip

\noindent {\sl Proof.} (i) We start with formulas from Theorem \eqref{thm.conv.3} where it is proved that
\begin{equation*}
|u(x,t)-M\,G_t(x)+\sum_i {\mathcal N}_{1,i}\,\partial_{x_i}G_t(x)|= D_2
\end{equation*}
with
$$
D_2= (4\pi t)^{-(N+2)/2} \int_{\ren} \int_0^1 |y|^2 u_0(y)\left(-\frac12 + (\langle \widehat y, \xi \rangle)^2 \right)  e^{-\xi^2}\,dyds,
$$
and $\xi=(x-sy)/2\sqrt{t}$. Multiplying by $|x|$ and integrating we have
$$
\int_{\ren}  D_2 \,|x|dx=(4\pi t)^{-(N+2)/2} \iiint  |y|^2 u_0(y)\left(-\frac12 + (\langle \widehat y, \xi \rangle)^2 \right)  e^{-\xi^2}\,|x|dx\,dyds.
$$
Writing now $|x|\le 2|\xi| \sqrt{t}+ s|y|$, we split the upper estimate of this integral into $I_1 + I_2$, where
$$
I_1=C t^{-(N+1)/2} \iiint  |y|^2 u_0(y)\left(-\frac12 + (\langle \widehat y, \xi \rangle)^2 \right)  e^{-\xi^2}\,2|\xi|dx\,dyds\,,
$$
and
$$
I_2= C t^{-(N+2)/2} \iiint  |y|^2 u_0(y)\left(-\frac12 + (\langle \widehat y, \xi \rangle)^2 \right)  e^{-\xi^2}\,s|y|dx\,dyds\,.
$$
 As for the first integral, the separate integrals are bounded and only the $\xi$ integral gets a time factor $t^{N/2}$ from  $dx=t^{N/2}d\xi$, so that
$$
I_1\le C\,t^{-1/2}.
$$
The second integral easily gives $I_2\le C\,t^{-1}$ by the assumption on the third moment of $u_0$. The proof is complete. \qed

\medskip

\noindent {\bf General conclusion}. All integrable solutions of HE in the whole space converge to Gaussian (in the renormalized forms we have written) if the initial mass is finite, $ u_0 \in L^1 (\ren) $. But the speed with which they do depends on how much initial mass is located far away, in other colloquial words for probabilists, on ``how  populated the tails are''. The quantitative versions we have established use mainly the moments of order 1 and 2. The moment of order 2 is called in Probability the (square of) the standard deviation. We remind the reader that not all probability distributions have a finite standard deviation (see Exercise 20).


\section{Functional Analysis approach for Heat Equations}\label{sec.conv.fa}

 We are going to use energy functions of different types to study the evolution of dissipation equations.  The basic equation is the classical heat equation, but the scope is quite general. Our aim is not to establish the convergence of general solutions to the fundamental solution (which is well done by other methods, as we have shown), but a bit more, namely, to find the speed  of convergence. After change of variables (renormalization) this reads as rate of convergence to equilibrium and relies on  functional inequalities. These functional inequalities also play an important in other areas.

The methods we will introduce next will apply to more general linear parabolic equations that generate semigroups. The method also works for equations evolving on manifolds as a base space. Since around the year 2000 we have been studying these questions for nonlinear diffusion equations. The main nonlinear models are: the porous medium equation,  the  fast diffusion equation, the $p$-Laplacian evolution equation, the chemotaxis system, some thin film equations, ...
Recently, the fractional heat equation and various fractional porous medium equations have been much studied.


\hiddensubsection{Heat Equation Transformations}

Take the classical Heat Equation posed in the whole space $\ren$ for $\tau>0$:
\begin{equation*}
u_\tau =\frac12 \Delta_y u
\end{equation*}
with notation $u=u(y,\tau)$ that is useful since we want to save the standard notation $(x,t)$ for later use. We know the (self-similar) fundamental solution, also called the evolution Gaussian profile
$$
U(y,\tau)=C\,\tau^{-N/2}e^{-y^2/2\tau}.
$$
It  was  proved in previous sections that this Gaussian is an attractor for all solutions in its {\sl basin of attraction}, consisting on all solutions with initial data that belong to $L^1(\ren)$ with integral $M=1$. See Sections  \ref{sec.scaling}, \ref{sec.conv1}, and \ref{sec.conv2}.

\noindent {\bf Remark.} We have inserted a harmless factor $1/2$ in front of the left-hand side following the probabilistic convention in order to get a Gaussian with clean exponent $-|y|^2/2t$ which has
standard deviation $t^{1/2}$ with no extra factors. Eliminating the prefactor leads to the exponential expression with exponent $-|y|^2/4t$, usual in PDE books. Our convention leads to some other simpler constants.

\medskip

\noindent $\bullet$ {\bf Fokker-Planck equation.} It is the first step in this approach to the asymptotic study. The scaling on the variables $u$ and $y$ to factor out the expected size of both which must mimic the Gaussian sizes, and then take logarithmic scale for the new time
$$
 u(y,\tau)=v(x,t)\,(1+\tau)^{-N/2}, \qquad 2t= \log(1+\tau).
$$
After some simple computations this leads to the well-known \sl Fokker-Plank equation \rm for $v(x,t)$:
\begin{equation}
v_t= \Delta_x v + \nabla_x\cdot(x\,v)=\nabla_x\left(\nabla_x v + xv\right)\,.
\end{equation}
We can write it as \ $v_t=L_1(v),$ where the Fokker-Planck operator $L_1=L_{FP}$ can be written in more explicit form as
$$
L_1(v)=\Delta v + x\cdot\nabla v + N\,v\,.
$$
We check now that when we look for stationary solutions by putting $v_t=0$ we
get as easiest case the  equation $\nabla v+ x\cdot v=0$ (after cancelling a divergence). Integrating it under the radial symmetry assumption is the simplest way to get the
Gaussian distribution $G=c\,e^{-x^2/2}$, and indirectly, the fundamental solution of the original heat  equation. We choose the constant $c=(2\pi)^{-N/2}$ to normalize $\int G\,dx=1$.

The asymptotic result we are aiming at consists precisely of proving that when $v_0(x)$ is integrable with mass 1 then $v(x,t)$ will tend to $G$ as $t\to\infty$. For a general presentation of the FP equation see \cite{Risken}. We will keep the notation $v(x,t)$ for the solutions of the Fokker-Planck equation throughout this section.

\medskip

\noindent $\bullet$ {\bf The Ornstein-Uhlenbeck semigroup.}
(i) In order to study relative error convergence it seems reasonable to pass to the quotient  $w = v/G$, where $G$ is the stationary state, to get the \sl Ornstein-Uhlenbeck \rm version
 \begin{equation}
w_t= L_2(w)= G^{-1}\,\nabla\cdot\big(\,G\,\nabla w\,\big)=\Delta w -  x\cdot\nabla w\,,
\end{equation}
a symmetrically weighted heat equation with Gaussian weight. Note that the corresponding stationary solution is now $W=1$, much easier.

\medskip

(ii) The two-term form of the r.h.s. looks easier, with a diffusion and a convection term.  Indeed, the weighted form of the {\sl Ornstein-Uhlenbeck operator} $L_2=L_{OU}$ is very convenient for  our calculations. To begin with, it allows to prove the symmetry of the operator in the weighted space $X=L^2(Gdx)$: for every two convenient functions $w_1$ and $w_2$ we have
\begin{equation}
\int_{\ren} (L_2w_1)\,w_2\,G dx=\int_{\ren} w_1\,(L_2w_2)\,Gdx=
-\int_{\ren} \langle \nabla w_1,\nabla w_2\rangle Gdx.
\end{equation}
It seems natural to introduce the Gaussian measure $d\mu=Gdx$ as a reference measure in the calculations. The important consequence of this computation is that $A=-L_2$ is a positive and self-adjoint operator in the Hilbert space $X=L^2(d\mu)$. This is a rather large space that includes all functions with polynomial growth. We will keep the notation $w(x,t)$ for the solutions of the Ornstein-Uhlenbeck equation throughout this section.

\medskip

(iii) We may also observe that $L_1(v)=G\,L_2(v/G)$ and that $L_2$ is the adjoint to $L_1$ in the sense that for conveniently smooth and decaying functions
$$
\int_{\ren} (L_1 v)w\, dx=\int_{\ren} v\,(L_2w)dx.
$$
More formally, we can consider the duality between the spaces $X=L^2(Gdx)$ and $X'=L^2(G^{-1}dx)$ given precisely by the integral of the product and the operators are adjoint.

\medskip

(iv) Finally, to complete the comparison we can write the Fokker-Planck equation as
$$
\partial_t v= \nabla \cdot (G\nabla (v/G))\,.
$$
The analogy says that this operator is negative and self-adjoint in the stranger space
$X_1=L^2(G^{-1}dx)$, that is much smaller than $L^2(\ren)$. For a detailed mathematical presentation of the Ornstein-Uhlenbeck semigroup we refer to \cite{AnSj}.

\medskip

\noindent $\bullet$ {\bf The Hamiltonian connection.} Start from the Fokker-Planck equation and use now the change of variables \ $v=zG^{1/2}$. Then,
$$
G^{1/2}z_t= \nabla\cdot(G\nabla (z G^{-1/2})).
$$
We have
$$
\nabla\cdot(G\nabla (z G^{-1/2}))=G^{1/2}\Delta z+ \nabla G^{1/2} \cdot\nabla z+G\nabla G^{-1/2}\cdot \nabla z+ z \nabla (G \nabla G^{-1/2}),
$$
 so that the equation for $z$ becomes:
$$
z_t= \Delta z- V(x)z, \qquad V(x)=G^{-1/2}\Delta G^{1/2}\,,
$$
that we may write as a  real Schr\"odinger Equation $z_t= L_3(z)=-H(z)$ with Hamiltonian operator
$$
H(z)=-\Delta z+V(x)z, \qquad V=\frac14 |x|^2-\frac{N}2,
$$
In calculating the Schr\"odinger potential we have used $G^{1/2}=e^{-|x|^2/4}$. Operator $H$ is directly symmetric in $L^2$ with no weight. The fact that it is positive is not clear from the formulas but it will follow from the equivalence with the Ornstein-Uhlenbeck operator.

\medskip

\noindent $\bullet$ The equivalence of this equation with the former ones comes from the transformation formulas
$$
L_1(v)=G^{1/2}L_3(v/G^{1/2})=G\,L_2(v/G),
$$
and also that if $v_i=G^{1/2}z_i=G w_i$, $i=1,2$, we get
$$
\int_{\ren} v_1\,v_2 \, G^{-1}dx=\int_{\ren} z_1\,z_2\,dx=\int_{\ren} w_1\,w_2\,Gdx,
$$
which allows to show that all three operators $L_i$ are self-adjoint and dissipative since we have
already proved it for $L_2$. The interesting remark for the Schr\"odinger representation is that it does not need any weighted space, $X_3=L^2(\ren)$.

\medskip

The rich equivalence among the three equations and also with the heat equation is a favorite topic in Linear Diffusion  and Semigroup Theory.

\medskip

\noindent $\bullet$ {\bf General Fokker-Planck Equation.}
It is based on generalising the coefficient $x$ of the convection term into a more general term that is the gradient of a potential that we call $S(x)$.\footnote{The standard notation is $U(x)$ but we will change the notation  here to $S(x)$ to avoid confusion with other uses of $U$.} The General Fokker-Planck equation (GFP) reads
\begin{equation}
v_t=\Delta v+ \nabla \cdot (\nabla S\,v)\,.
\end{equation}
Standard assumption is that $S$ must be a positive and convex function in $\ren$, called the potential. The stationary state is now $\widetilde G(x)=C\,e^{-S(x)}$, and the equation reads then $v_t=L_1(v)$
with
$$
L_1(v)=\nabla \cdot (\widetilde G\,\nabla (v/\widetilde G)).
$$
The other equations, General Ornstein-Uhlenbeck and General Hamiltonian Equation, follow in the same way as before using only the expressions in terms of $G$ that is replaced by $\widetilde G$. The weighted scalar products have no difference and the relation of norms still holds. In the Hamiltonian representation we get a potential (put $G=e^{-S}$)
$$
V(x)=-G^{-1/2}\Delta (G^{1/2})= \frac14|\nabla S|^2-\frac12\Delta S\,.
$$

On the other hand, the analysis of the complete spectrum is not possible unless we have very particular cases of potentials $S$. Moreover, the connection with a renormalization of the heat equation is completely lost.


 \hiddensubsection{Asymptotic Energy Method via the Ornstein-Uhlenbeck Equation}

The Ornstein-Uhlenbeck formulation allows for a very clear and simple treatment of the problem of convergence with rate to the Gaussian profile. We may assume without lack of generality that
$$
\int w\,d\mu=\int v\,dx=\int u\,dy=1\,,
 $$
 with the above notations for $u, v,$ and $w$.  We now make a simple but crucial calculation on the time decay of the energy for the OUE:
\begin{equation}
{\mathcal F}(w(t))=\int_{\ren} |w-1|^2\,G\,dx, \quad \frac{d{\cal F}(w(t))}{dt} =-2\int_{\ren}|\nabla w|^2\,G\,dx
=-{\mathcal D}(w(t)).
\end{equation}
We can now use a result from abstract functional analysis: the {\bf Gaussian Poincar\'e inequality } with measure  $d\mu=G(x)\,dx$:
$$
\int_{\ren} w^2d\mu- \left(\int_{\ren}w\,d\mu\right)^2 \le C_{GP} \int_{\ren} |\nabla w|^2\,d\mu
$$
The sharp constant in this inequality is precisely $C_{GP}=1,$ with no dependence on dimension. Moreover the functions that realize the optimal constant are $w(x)=x_i$ for any $i=1,2,\dots,N$. \footnote{This is an old inequality in the folklore of Hermite polynomials, and probably was known in one dimension to both mathematicians and physicists in the 1930's in relation to eigenvalue problems, as mentioned in \cite{Beck}.} We will give below a proof based on the analysis of the spectrum of the Ornstein-Uhlenbeck operator.

Then, the left-hand side is just ${\cal F}$ and the inequality implies
$$
-\frac{d}{dt}\,{\mathcal F}(w(t))\ge 2{\mathcal F}(w(t)),
$$
which after integration gives ${\mathcal F}(w(t))\le {\mathcal F}(w_0)e^{-2t}$, i.e.:
\begin{equation*}
\int_{\ren}|w-1|^2\,\rd\mu \le {\rm e}^{-2t}\,
\int_{\ren}|w_0-1|^2\,\rd\mu\quad\forall\; t\ge 0
\end{equation*}

We have proved the convergence to equilibrium in the following form.

\begin{theorem} Under the assumptions on the initial data
$$
\int_{\ren}w_0\,\rd\mu =1, \qquad \int_{\ren}w_0^2\,\rd\mu <\infty\,,
$$
the solutions of the OUE satisfy the following stabilization estimate
$$
\|w(t)-1\|_{L^2(d\mu)}\le \|w_0-1\|_{L^2(d\mu)}\,{\rm e}^{-t}.
$$
\end{theorem}

This is the first of the well known Functional  Estimates for the solutions to the HE.
The renormalization  $\int w_0\, d\mu=1$ is no restriction.

In terms of $ v $ the hypotheses are $ \int v_0 \,dx = 1 $ and $ \| v_0-G \|_{L^2 (G^{-1}dx)} <\infty $. The weight is now $K=G^{-1}$  and the measure $ d \nu = K \,dx $, which behaves inversely in infinity. The result is
$$
\int_{\ren}|v-G|^2\,d\nu \le {\rm e}^{-2t}\,
\int_{\re^d}|v_0-G|^2d\nu \quad\forall\; t\ge 0.
$$

Recall that we have a new logarithmic time $t=\log(1+\tau).$ The rate of convergence in real time in both variables is then $O((1+\tau)^{-1/2})=O(\tau^{-1/2})$ as $\tau\to \infty$.

\medskip

\noindent $\bullet$ Then by regularity theory (regularizing effect from $L^1$ to $L^\infty$) for the heat equation and the other equations,  we get estimates in the sup norm with similar relative rates at least locally in space. But  note that since $w=xe^{-t}$ is a solution of the OUE,  we cannot get a uniform estimate in $L^\infty$ of the  Ornstein-Uhlenbeck variable.


\medskip

\section{Calculation of spectrum. Refined asymptotics}\label{sec.spec}

We will proceed with a further step in the analysis to get more precise asymptotics. Indeed, the knowledge of the spectrum of the equivalent operators allows to obtain a complete description of the long-time behaviour in weighted spaces. This is a well-known fact in the study of the HE posed in bounded domains, that has a parallel here.

\hiddensubsection{Spectrum} We will make all computations on the Ornstein-Uhlenbeck operator. The Gaussian Poincar\'e inequality is a simple consequence of the following analysis of the spectrum of the Ornstein-Uhlenbeck Operator in $L^2(d\mu)$:

(a) Since the FP and the OU operators  have a compact inverse, we conclude that they have a discrete spectrum, \cite{Kavian}. The ground state of the OUE is formed by the constant function $w=1$ (which comes from the Gaussian function $v=G$ for the FPE) with eigenvalue $\lambda_0=0$. The next eigenfunctions are the coordinate functions $\phi_i=x_i$ corresponding to the  eigenvalue $\lambda_1=1$ with multiplicity $N$.

(b) In 1D we find the rest of the eigenfunctions and eigenvectors as the family of Hermite polynomials, given by the compact formula
\begin{equation}
H_k(x)= (-1)^k G(x)^{-1} (d/dx)^k G(x)\,,
\end{equation}
that tells much about how we will see them. Indeed, the formula can be derived from the fact that the derivatives in $y$ of the Gaussian evolution solution $U(y,\tau)$ are still solutions of the heat equation with different decay rate. Passing to the FPE we conclude that the derivatives $D^k G$ are eigenfunctions of $L_1$ (here $D=d/dx$). It is easy to see that these solutions  have the form $H_k(x)G(x)$ where $H_k$ is a polynomial of degree $k$ (proof by induction). Passing to the OUE we get the formula above. More precisely, we have the recursion formula:
$$
G H_{k+1}=-d/dx(H_k G)=-(H_k'-xH_k) G, \quad H_{k+1}=(x-\frac{d}{dx})H_k\,.
$$
The first members of the family $H_k$ are $1, x, x^2-1, \dots$

The corresponding eigenvalue to $H_k$ for $L_2=L_{OU}$ is $\lambda_k=k$. This can be seen from the heat equation formula since differentiating in $x$ adds a factor $\tau^{-1/2}$ to the decay, which goes over as $e^{-t}$ for every derivative we take. Induction proof at the FP level: if we assume that
$$
L_1(D_kG)=(G(D_kG/G)')'=(D_kG)''+x(D_kG)'=-\lambda_kD_kG\,,
$$
then
$$
L_1(D_{k+1}G)=(D_kG)'''+x(D_kG)''=-\lambda_kD_kG'-(x(D_kG)')'+ x(D_kG)''=-(\lambda_k+1)D_{k+1}G.
$$
Therefore, $\lambda_{k+1}=\lambda_k+1$. Since $\lambda_0=0$ the proof is done.

It is then proved that the $H_k(x)$ form a basis in $L^2(d\mu)$ in 1D.

(c) For several dimensions we have the functorial property: if $x=(x_1,x_2) $ and $w(x)=w_1(x_1)w_2(x_2)$ we get
$$
L_2(w)=w_1\,L_2(w_2)+ L(w_1)\,w_2\,.
$$
This produces new eigenfunctions in higher dimensions, and it gives for the product function the sum of the eigenvalues: $\lambda(w)=\lambda(w_1)+\lambda(w_2)$. The set of combinations generates a base of eigenfunctions, this is essentially due to Fubini's theorem, cf. \cite{Bell}. Our account is very short but we consider this part an extension, and it is well documented in the corresponding literature.

\medskip


 \hiddensubsection{Refined asymptotics}

 From the spectrum we can get a very precise description of the convergence in the weighted spaces
by using the equivalent of the Fourier analysis on bounded domains. The meaning of the coefficients for the original equation has to be understood.

We get
$$
w(x,t)=\sum_\alpha c_\alpha\, H_\alpha(x)e^{-kt}, \quad v(x,t)=\sum_\alpha c_\alpha \,\partial_\alpha G(x)e^{-kt}
$$
where $\alpha$ is an $N$-multi-index, $k=|\alpha|$, and $H_\alpha$ is the corresponding multidimensional Hermite polynomial, after renormalization in $L^2(d\mu)$. We have
$$
c_\alpha=\frac{\langle w_0, H_\alpha\rangle_{L^2_\mu}}{\langle H_\alpha\, H_\alpha\rangle_{L^2_\mu}}=\frac{\int v_0(x)H_\alpha(x)\,d\mu(x)}{\int H^2_\alpha(x)\,d\mu(x)},
$$
\normalcolor
so that $c_0$ is the  mass and for $k=1$, and $c_i$ are the first coordinate moments $\int v_0(x)x_i\,dx$ after normalization, and so on.
The convergence of the $w$ series holds in $L^2(d\mu)$, with errors of the order of the first term that is left out. The convergence of the $v$ series holds in $L^2(G^{-1}dx)$, with errors of the order of the first  term that is left out.

\section{Convergence via the Boltzmann entropy approach}\label{sec.entropy}

  There is another approach for the convergence to the Gaussian that starts the analysis from Boltzmann's ideas on entropy dissipation. We start now from the Fokker-Planck equation $v_t=\Delta v+ \nabla\cdot(xv)$ and consider the functional called  {\bf entropy}
$$
{\cal E}(v):=\int_{\ren} v\,\log(v/G)\,dx=\int_{\ren} v\,\log(v)dx+\frac12\int_{\ren}  x^2 v\,dx +C\,.
$$
and we assume that the data are such that the initial entropy is finite.  We recall that no decay is possible without some restriction on the data.

Differentiating along the flow (i.e., for a solution) leads to
$$
\frac{d{\cal E}(v)}{dt}=-{\cal I}(v), \quad
{\cal I}(u)=\int_{\ren} v\,\left|\frac{\nabla v}v+ x \right|^2\,dx=\int_{\ren} v\,|\nabla \log(v/G)|^2\,dx\,.
$$
For some reasons the dissipation ${\cal I}(v)$ is called Fisher information. Let us continue the proof. Putting now $v=Gf^2$ we find that
$$
{\cal E}(v)=2\int_{\ren} f^2\,\log(f)\,d\mu, \quad {\cal I}(v)=4\int_{\ren} |\nabla f|^2\,d\mu.
$$
The famous {\sl logarithmic Sobolev inequality} proved by Gross in 1975, \cite{Gross75}, says than that (for all suitable functions, not only solutions)
$$
{\cal E}\le \frac 12 {\cal I}\,,
$$
and we obtain the decay ${\cal E}(t)\le {\cal E}(0)\,e^{-2t}.$ This means a precise decay for the entropy functional. The calculations are justified for smooth solutions, and then we can pass to the limit for general solutions with finite entropy.

In order to obtain decay in standard norms, there are formulas connecting the entropy with other norms, like the  Cziszar-Kullback inequality that implies that
that
$$
\|f-G\|^2_{L^1(\ren)}\le 2{\cal E}(f, G),
$$
for any $f, G\in L^1(\ren)$ positive with equal mass, see \cite{Barron}. This is paper is a very good early reference to the subject of entropies and the central limit theorem.

There are many works dealing with the use of functionals and functional inequalities to arrive at asymptotic behaviour results plus a rate of convergence for this kind of equations. Let us mention here \cite{A4, BE84, Chafai, MV, Tosc} and the references to be mentioned in Section \ref{sec.appl2}.


 \hiddensubsection{About entropy in Physics}

Entropy has been introduced as a state function in thermodynamics
by R. Clausius in 1865, in the framework of the second law of thermodynamics, in order to interpret the results of S. Carnot.

A statistical physics approach: Boltzmann's formula (1877) defines the entropy of a physical system in terms of a counting of its micro-states. Boltzmann's equation:
$$\partial_t f + v \cdot\nabla_x f = Q(f, f)\,.
$$
 It describes the evolution of a gas of particles having binary collisions at
the kinetic level; $f(t,x,v)$ is a time dependent distribution function
(probability density) defined on the phase space $(x,v)\in \ren\times\ren$. The Boltzmann entropy: \ $H[f] :=\iint f \log (f) dx dv$  measures irreversibility. The famous H-Theorem (1872) says that
$$
\frac{d}{dt}H[f] =\iint Q(f, f) log (f) dx dv\le 0\,.
$$
Other approaches to thermodynamic entropy are due to Carath\'eodory (1908), Lieb-Yngvason  (1997),... see \cite{LY}.

An important version of entropy appears in Information Theory. In 1948, while working at Bell Telephone Laboratories Claude Shannon, an electrical engineer, set out to mathematically quantify the statistical nature of ``lost information'' in phone-line signals (cf. Wikipedia article). He arrived at an analog to thermodynamic entropy for use in information entropy.

There is also a concept of entropy in probability theory (with reference to an arbitrary measure).

\section{Brief review of other heat equation problems}\label{sec.appl1}

The methods presented above have been applied to prove convergence to a distinguished solution (that plays the role of the Gaussian fundamental solution) in different contexts. We have not mentioned some other methods like the transport method of Jordan-Kinderlehrer-Otto \cite{JKO98}, 1998, where the Fokker-Planck equation is interpreted as the steepest descent for a free energy related to Boltzmann-Gibbs entropy, taken with respect to the Wasserstein metric. This novel technique, based on mass transportation \cite{VillaniOT}, has played an increasing role since then

 \medskip

\hiddensubsection{Equation with forcing} A modification that still keeps the flavor of this presentation consists of considering a forcing term
\begin{equation}
\partial_t u=\Delta u +f\,,
\end{equation}
where $f$ is an integrable function of $(x,t)\in Q_T=\ren\times(0,T)$, $T>0$. It can be easily proved that a representation formula holds, \cite{EvansPDE}. From it we can derive asymptotic results that we leave as exercises.

\medskip

\noindent {\bf Exercise 21.} (i) Prove that the  $L^1$ estimate of Theorem \ref{main.convthm}, formula \eqref{cr.l1},  holds  if $f\in L^1(Q_\infty)$, and we take as $M$ the accumulated mass  defined as
\begin{equation}
M=\int_{\ren} u_0(x)\,dx+\iint_{Q_\infty} f(x,t)\,dxdt.
\end{equation}

\noindent (ii) The $L^\infty$ statement \eqref{cr.linf} needs some further decay condition on $f$ as $t\to\infty$. Prove that it holds e.g. if $\|f(\cdot,t)\|_\infty \le Ct^{-\gamma}$ with $\gamma\ge 1+N/2$.

\medskip

 The results of Section \ref{sec.conv1} can be repeated under conditions that we invite the reader to provide.

\medskip

\noindent $\bullet$ Let us continue with a variation of the heat equation with forcing. In fact, there are many studies where the forcing term takes the form $f=f(u)$, and they fall into what is called reaction-diffusion,
\cite{Smo82}. The simplest case correspond to linear forcing, $f(u)=\kappa u$ with $\kappa\ne 0$. The modifications  in the analysis are minimal since the change of variables $v(x,t)=u(x,t)\,e^{-\kappa t}$ transforms equation
\begin{equation}
\partial_t u=\Delta u +\kappa u\,,
\end{equation}
into the classical form $v_t=\Delta v$, to which our results apply. We thus conclude that
for very large $t$ we have the asymptotic behaviour
\begin{equation}
u(x,t)\sim M\,G(x,t)\,e^{\kappa t}\,,
\end{equation}
which preserves the Gaussian profile as asymptotic shape but not the decay rates in time. The reader is asked to write the precise theorems using the results of Sections \ref{sec.scaling}, \ref{sec.conv1}, and \ref{sec.conv2}.

 \medskip

\hiddensubsection{Dipoles and related issues}
We have already seen that the heat equation with signed $L^1$ initial data and zero mass, i.\,e., $\int u_0(x)\,dx=0$ in the sense that the term representing the Gaussian approximation vanishes, so that the
rate of decay as $t\to\infty$ is faster and the first approximation is given by a first term that combines partial derivatives of the Gaussian :
\begin{equation}
u(x,t)\sim D_{\bf v} G_t(x)\,, \qquad {\bf v}=\sum_i {\cal N}_{1,i}(u_0)\,e_i\,,
\end{equation}
of course  under the condition that this vector $\bf v$ does not vanish. The $e_i$ are the canonical basis and $D_{\bf v}$ denotes directional derivative. Therefore, $u(x,t)= O(t^{-(n+1)/2})$.

\medskip

\noindent $\bullet$ There is a very interesting application of this result in $N=1$. Indeed, we may  solve the problem of asymptotic behaviour of the solutions of the heat equation posed in a half line $I=(0,\infty)\subset \re$ for $t>0$ with lateral Dirichlet data $u(0,t)=0$, let us call it (DP-HE-HL).
We assume that $u_0\in L^1(I)$ and $u_0\ge 0 $ (the last assumption is made for simplicity). The idea is tho extend the initial data to the whole line by putting
$$
\ou_0(x)=u_0(x) \quad \mbox{if } \ x>0, \qquad
\ou_0(x)=-u_0(-x) \quad \mbox{if } \ x<0
$$
(called anti-symmetric reflection). Note that $\ou\in L^1(\re)$,  its total mass is zero, and the first moment $2{\cal N}_1$ is not zero. Solving the heat equation with the usual representation formula we obtain a solution $\ou(x,t)$ defined
for $x\in\re$ and $t>0$. This solution must be antisymmetric, $\ou(x,t)=-\ou(-x,t)$ by the form of the data and an elementary symmetry property of the heat equation. Restricting $\ou$ to $x>0$ we find a unique solution of Problem  (DP-HE-HL). We can now copy our asymptotic results to $\ou$ and translate them to $u$. Let us write
\begin{equation}\label{dipole}
D(x,t)=-\partial_{x}G_t(x)=\frac{x}{t^{3/2}}e^{-|x|^2/4t}\,.
\end{equation}
Theorems \ref{thm.conv.3} and \ref{thm.conv.4} imply that $D(x,t)$ is the asymptotic attractor of the evolution in the half line.

\begin{theorem}\label{thm.dipoleconv} Let us assume that $u_0\in L^1(I)$  and the first moment in $I$
\begin{equation}
{\mathcal N}_{1}=\int_0^\infty y\,u_0(y)\,dy
\end{equation}
is finite, as well as  the second moment :
${\mathcal N}_2=\int_0^\infty |y^2\,u_0|(y)\,dy<\infty.$
Then, we  get the following convergence formulas with rate for the solutions of (DP-HE-HL)
\begin{eqnarray*}\label{conv.dip}
\displaystyle \|u(x,t) - 2{\mathcal N}_{1}\, D(x,t)\|_{L^1(I)}\le C\,{\mathcal N}_2 \,t^{-1}\,,\\
\displaystyle t^{1/2}|u(x,t)- 2{\mathcal N}_{1}\, D(x,t)|\le
C{\mathcal N}_2\,t^{-1}\,.
\end{eqnarray*}
The rate $O(t^{-1})$ is optimal under such assumptions. If the third moment is finite we also have
\begin{equation}
\displaystyle \|u(x,t) - 2{\mathcal N}_{1}\, D(x,t)\|_{L^1(I; |x|dx)}\le C\,({\mathcal N}_2 + {\mathcal N}_3)\,t^{-1/2}.
\end{equation}
\end{theorem}

Function $D(x,t)$ given by \eqref{dipole} is called the {\bf dipole solution} because it takes the derivative of the unit Dirac delta as initial data. It has a constant-in-time first moment that characterizes its strength. Not that the mass in $I$, $M_+=\int_0^\infty u(x,t)\,dx$, is not conserved in time, but decays like $O(t^{-1/2})$.
\footnote{Dipole solutions appear often in Physics, specially in Electromagnetism.}

We can derive a more general convergence result

\begin{theorem}\label{thm.dipoleconv2} Let us assume that $u_0\in L^1(I)$  and the first moment in $I$,  ${\mathcal N}_{1}$, is finite.
Then,
\begin{equation}
\displaystyle \|u(x,t) - 2{\mathcal N}_{1}\, D(x,t)\|_{L^1(I; |x|dx)}\to 0
\end{equation}
as $t\to \infty$.
\end{theorem}

\noindent {\sl Proof.} (i)  Assume that $u\ge 0$. Approximate $u_0$ by a compactly supported function $\widetilde u_0\le u_0$, so that the error in the first moment is $\ve$.  Then,
$$
\|u(x,t) - \widetilde u(x,t)\|_{L^1(I; |x|dx)}=\int_0^\infty (u(x,t)-\widetilde u(x,t)) xdx \le \ve
$$
because this expression is  conserved with time for solutions of  Problem  (DP-HE-HL) in $I$ (check this!).
The previous theorem for the solution $\widetilde u$ gives
$$
\|\widetilde u(x,t) - 2{\widetilde {\mathcal N}}_{1}\, D(x,t)\|_{L^1(I; |x|dx)}\le K\,t^{-1/2}.
$$
that can be made less than $\ve$ if $t$ is large enough. Finally, $| {\mathcal N}_{1}-{\widetilde {\mathcal N}}_{1}| \le \ve$ and $D$ has a finite moment. Combining all this, the result holds.

(ii) For signed data, split into positive and negative part, and then combine the results. \qed

\medskip

\noindent $\bullet$ Similar formulas hold for the heat equation posed in a half space $\Omega$ in $N$ dimensions. After rotation and translation we may take \ $\Omega=\{x\in \ren: x_1>0\}$. We consider zero lateral Dirichlet data: $u(x,t)=0$ for $x=(0,x_2,\dots,x_N)$. n this case a multi-dimensional dipole solution appears, $D=-\partial_{x_1}G_t$.
Problems in half spaces can be solved too. The equivalent of Theorem \ref{thm.dipoleconv} holds. We leave the easy details to the reader.

Finally, solutions of the type $D_{12}=\partial_{x_1}\partial_{x_2}G_t$ are attractors for problems posed in quadrant domains $\{x\in \ren: x_1>0, x_2>0\}$. And so on.

\medskip

\noindent $\bullet$ By using the  symmetric extension instead of the anti-symmetric one we can solve the heat equation posed in a half line $I=(0,\infty)\subset \re$ for $t>0$ with lateral Neumann data $u_x(0,t)=0$, let us call it (NP-HE-HL). Let the reader fill in the details if needed.

\noindent $\bullet$ A related problem in several dimensions occurs when the domain is an exterior domain with one or several holes and appropriate boundary conditions. Thus, with zero Dirichlet boundary conditions and integrable data convergence to the Gaussian holds,  while in 1D we fall back into the dipole problem.

\hiddensubsection{Other problems in subsets of $\re^N$}
There are a number of problems involving the heat equation that have been studied in great detail, like the heat equation posed in a bounded domain with boundary conditions of different types (Dirichlet, Neumann, mixed, or other), but these settings lead to quite different results that depart too much from the picture presented here, so we will not comment on them.

There are also equations with coefficients or weights;  they form a large topic that leads also very far from the present presentation.

\hiddensubsection{Heat equation on manifolds}
 The construction of the heat equation has been carefully studied when the equation is posed on a Riemannian manifold $(M^N,g)$. The equation then takes the form
  \begin{equation}
\partial_t u(x,t)= \Delta_{g}\, u(x,t)={|g(x)|^{-N/2}}\sum_{i,j=1}^{N}\partial_i \left(g(x)^{ij}|g(x)|^{N/2}\partial_j u(x,t)\right)\,,
\end{equation}
where $g_{ij}$ is the metric tensor, $g^{ij}$ its inverse, $|g|$  its determinant, so that $\Delta_{g}$ is the Laplace-Beltrami operator, \cite{Gr2009}. Two particular manifolds are specially relevant because their internal symmetries  and homogeneity make the theory specially strong and mathematically appealing: the $N$-dimensional sphere $\mathbb{S}^N$ and the hyperbolic space $\mathbb{H}^N$. The heat flow on the former is easily shown to stabilize to a constant (much like a Neumann problem in $\re^N$). The flow on  $\mathbb{H}^N$ is more interesting, and typical solutions with finite mass converge to a modified Gaussian function, the hyperbolic fundamental solution, that is described in detail in \cite{Grig98}, see also \cite{Grig87}.

\section{Application to other diffusion equations}\label{sec.appl2}

We now examine some nonlinear diffusion equations where similar methods and results have been successfully proved in the last half century, where the Gaussian profile is replaced by some other attractive object.

\medskip

\noindent $\bullet$ A prominent example that has been much studied is the Porous Medium Equation, $\partial_t u= \Delta u^m$, or better $\partial_t u= \Delta (|u|^{m-1}u)$, $m>1$. The asymptotic study depends crucially on the existence and properties of a  distinguished family of solutions, the Barenblatt solutions, \cite{Bar52}, that are compactly supported and self-similar, one for every mass $M$.
Explicit formulas exist for them (1950, 1952):
\begin{equation}\label{Barensol}
{\bf B}(x,t;M)= t^{-\alpha} {\bf F}(x/t^{\beta}), \quad {\bf
F}(\xi)=\left(C - k \xi^{2}\right)_+^{1/(m-1)}
\end{equation}
where $C$ is a free constant (to be determined by the mass $M>0$) and $k>0$ is a function of $m,N$ ($k=(m-1)/2m(N(m-1)+2)$). They replace the Gaussian fundamental solutions in the statement of the asymptotic theorems. Convergence of finite mass solutions to the Barenblatt solution with the same mass is proved by the scaling method in \cite{Vaz03}. A much earlier proof used a method of optimal upper bounds, \cite{FrKa80}, 1980.

The faster convergence of solutions of the porous medium equation with changing sign is studied in \cite{KV91}.  Here dipole solutions appear when the total mass is zero. Dipole solutions for nonlinear parabolic problems were also studied in \cite{BaZe57, GPV95, HuVa93}. Convergence to the Barenblatt kernel for the PME posed in $\ren$ minus one or several holes was studied in \cite{BQV2007, GG2007}.

These methods did not produce rates of convergence, but such rates were established using the entropy method by Carrillo-Toscani in \cite{CaTo00}, 2000, via Bakry-Emery inequalities \cite{BE84, BGL}, and then in Del Pino-Dolbeault \cite{DPD-2002}, 2002, using Gagliardo-Nirenberg inequalities. Convergence in the sense of Wasserstein distances was introduced by F. Otto in 2001, \cite{Otto}. Entropy methods suitable for weighted porous media equations are used in \cite{DGGW}.

\medskip

\noindent $\bullet$ The Fast diffusion equation (i.\,e., the PME for $m<1$) behaves much like the porous medium equation for $m$ close to 1, even if the shape of the selfsimilar profile is different,
\begin{equation}\label{Baren.fast}
{\bf F}(\xi)=\left(C + k(m,N) \xi^{2}\right)^{-1/(1-m)}
\end{equation}
with fat tails as $|x|\to\infty$. This is not the case for $m<(N-2)/N$ due to the phenomenon of extinction in finite time. A detailed analysis of convergence to so-called pseudo-Barenblatt profiles is done in \cite{BBDGV2, BDGVProc} by an entropy method which  relies on some Hardy-Poincar\'e inequalities. This is anyway a quite different scenario.

\medskip

\noindent $\bullet$ The asymptotic convergence for the $p$-Laplacian equation $\partial_t u= \Delta_p u$ has been treated in \cite{KVplap} by the scaling method, after settling the uniqueness of the fundamental solution of Barenblatt type. This was done for $p>2$, and it extends to some $p<2$, to be precise to $2N/(N+1)<p<1$. The entropy method is used in \cite{DPD-2002cras, DPD-2003} using the analogy with the porous medium equation. The doubly nonlinear equation $\partial_t u= \Delta_p (u^m)$ was studied in \cite{Agueh} and subsequent works.

%

\medskip

\noindent $\bullet$ There many variants of heat equation with lower-order terms, either first order or zero order. If these terms are strong enough they will destroy the convergence towards some attractive solution with a Gaussian shape. This is a huge research field and we will give only some ideas. Maybe the best known models correspond to the case that can be written as $\partial_t u=\Delta u + f(u)$ with $f(u)=-u^p$,  $u\ge 0$; we then have so-called diffusion-absorption equations. Let us take and $p>1$, since the case $p=1$ was explained before.  Even if conservation of mass does not hold, convergence of finite-mass solutions to a Gaussian profile with some positive mass $M_\infty$ is proved when $p>p_c=(2+N)/N$. The limit case $p_c$ is very interesting and was studied in \cite{GmVeron}. For $p<p_c$ we enter into completely new asymptotic profiles. See also \cite{KPV89, VazIMA}.

There are many extensions of these ideas. The reaction cases, $f(u)=+u^p$, lead to the existence of blow-up in finite time, a huge topic that falls complete out of the scope of these notes, see \cite{GaVaz20014, QuiSou2007, SGKM}. Another reaction case that has attracted the attention of researchers is $f(u)=u(1-u)$, called the Fisher-KPP model, of interest in biology and chemistry; here the long term behaviour takes the form of expanding travelling waves and no trace of a Gaussian is seen,   \cite{K-P-P:art}.

\medskip

\noindent $\bullet$ Let us now examine some of the recent work on diffusion with fractional Laplace operators. The linear fractional heat equation $\partial_t u+ (-\Delta)^s u=0$ has a rather complete theory in the paper \cite{BSiV2016}. It is well known that  the self-similar fundamental solution exists for every $N\ge 1$, $0<s<1$ and has the form $P(x,t;s)=t^{-N/2s}F_s(|x|t^{-1/2s})$ where $F_s$ is a smooth and positive profile function with a fat tail as $|x|\to\infty$
\begin{equation}
F_s(\xi) \sim c(N,s)|\xi|^{-(N+2s)}\,,
\end{equation}
see \cite{BG1960}. Convergence to the self-similar fundamental solution can be proved for finite-mass solutions by the scaling method (no rates), or by the representation analysis (with rates). See  separate notes by the author. For the entropy method see \cite{BK2003, GentImb}.

\medskip

\noindent $\bullet$ We continue with nonlinear fractional heat equations of porous medium type. The model studied by Caffarelli and V\'azquez \cite{CV1} admits self-similar solutions that we may call fractional Barenblatt solutions \cite{BKM10, CV2, BIK2015}. The entropy method is used in \cite{CV2} to establish asymptotic convergence without rates. Rates in 1D were obtained in \cite{CH2013}. Convergence with rates in several dimensions is not known.

The alternative model of fractional porous equation, $\partial_t u+ (-\Delta)^s(|u|^{m-1}u)=0$ was studied in \cite{pqrv1, pqrv2}.  Unique fundamental solutions of Barenblatt type were described in \cite{VazBar2012}, where convergence (without rates) was proved by the scaling method.

\medskip

\noindent $\bullet$ There are a number of other equations that have been studied, like thin film equations \cite{CaTo01}, the Barenblatt equation of elastoplastic filtration {KPV91}, inhomogeneous heat or porous medium equations with weights, \cite{KRV10}, or chemotaxis models, like \cite{BDEF}. A very important topic is the study of the heat equation and the nonlinear diffusion models on manifolds, like the hyperbolic space, see \cite{Grig98, Vazhyp15}.

\medskip

\noindent $\bullet$   See \cite{VazCIME} for a  general presentation of linear and nonlinear diffusion equations including a detailed survey of recent research work.

\medskip

\section{Historical comments}

We add some  historical notes on the origins and development of the Gaussian function,  borrowed from Wikipedia and other widely available sources, with no claim to be a rigorous historical presentation. It seems that the so-called Gaussian function has its origin in Statistics. The 18th century statistician Abraham de Moivre, a Frenchman exiled in England, seems to have been the first person who noticed the existence of a bell-shaped curve as the limit of the probability distributions of repeated random trials done independently. He was led by the practical problems of calculating odds in gambling, not a every elevated motivation indeed. But he was a very fine mathematician, appreciated by Newton. The curve he discovered is now called the ``normal curve''.

\smallskip

One of the first applications of the normal distribution was to the analysis of errors of measurements made in astronomical observations. A century later than De Moivre, the mathematicians Adrain in 1808 and Gauss in 1809 developed independently the formula for the normal distribution and showed that errors were fit well by this distribution. The brilliant Gauss received much credit, Adrain's work was not known for many years.

\smallskip

This same distribution had been discovered by Laplace in 1778 when he derived the extremely important central limit theorem, a main topic of  this paper. Laplace showed that even if a distribution is not normally distributed, the means of repeated samples from the distribution would be very nearly normally distributed, and that the larger the sample size, the closer the distribution of means would be to a normal distribution.

\smallskip

The distribution appeared later in another disguise in Statistical Mechanics as the  Maxwell-Boltzmann distribution, shortly called the Maxwellian. The original derivation in 1860 by James Clerk Maxwell was an argument based on molecular collisions of the kinetic theory of gases as well as certain symmetries in the speed distribution function; Maxwell also gave an early argument that these molecular collisions entail a tendency towards equilibrium. After Maxwell, Ludwig Boltzmann in 1872 also derived the distribution on mechanical grounds and argued that gases should over time tend toward this distribution, due to collisions (see H-theorem).

\smallskip

The normal distribution has a wide implication to social issues. Thus, Qu\'etelet seems to have been the first to apply the normal distribution to human characteristics. He noted that characteristics such as height, weight, and strength were normally distributed.

\smallskip

Evidence for the Gaussian function as the fundamental solution of the heat equation came after the work of probabilists in the 20th century to establish the link between heat equation and Brownian diffusion, which is turn is the limit of discrete processes based on iterated random trials. The close connection between stochastic differential equations and parabolic partial differential equations is very much influenced by the role of the Gaussian function in both theories.

\

\

\noindent {\sc Acknowledgment.} Work partially supported by Spanish Project MTM2014-52240-P. These notes developed from Ph. D. courses and lectures given by the author at different events in recent years, the last one was the Annual Meeting of the Red de An\'alisis Funcional y Aplicaciones, held in C\'aceres, Spain, in March 2017.


\

\bibliographystyle{amsplain} 

    \


\

\noindent {\sc Address:}

\noindent Juan Luis V\'azquez. Departamento de Matem\'{a}ticas, Universidad
Aut\'{o}noma de Madrid,\\ Campus de Cantoblanco, 28049 Madrid, Spain.  e-mail address:~\texttt{juanluis.vazquez@uam.es}

\

\

\end{document}